\documentclass[12pt]{article}

\newtheorem{Df}{Definition}[section]
\newtheorem{Te}[Df]{Theorem}
\newtheorem{Po}[Df]{Proposition}
\newtheorem{Cr}[Df]{Corollary}
\newtheorem{Lm}[Df]{Lemma}
\newtheorem{Ca}[Df]{Claim}
\newtheorem{Cn}[Df]{Conjecture}

\newcommand{\Bdf}{\begin{Df}}
\newcommand{\Edf}{\end{Df}}
\newcommand{\Bte}{\begin{Te}}
\newcommand{\Ete}{\end{Te}}
\newcommand{\Bpo}{\begin{Po}}
\newcommand{\Epo}{\end{Po}}
\newcommand{\Bcr}{\begin{Cr}}
\newcommand{\Ecr}{\end{Cr}}
\newcommand{\Blm}{\begin{Lm}}
\newcommand{\Elm}{\end{Lm}}
\newcommand{\Bca}{\begin{Ca}}
\newcommand{\Eca}{\end{Ca}}
\newcommand{\Bcn}{\begin{Cn}}
\newcommand{\Ecn}{\end{Cn}}

\begin{document}

\title{\bf{Koszulity for nonquadratic algebras II}}
\author{
  \bf{Roland BERGER}\\
  LARAL, Facult\'e des Sciences et Techniques,\\
  23 rue P. Michelon,\\
  F-42023 Saint-Etienne Cedex 2, France\\
Roland.Berger@univ-st-etienne.fr\\
  }
\date{}
\maketitle

\begin{abstract}
It has been shown recently, in a joint work  with Michel
Dubois-Violette and Marc Wambst~\cite{bdvw:homog} (see also 
math.QA/0203035), that Koszul property of
$N$-homogeneous algebras (as defined in~\cite{rb:nonquad}) becomes
natural in a $N$-complex setting. A basic question is to define the differential of the bimodule Koszul complex of an $N$-homogeneous
algebra, e.g., for computing its Hochschild homology. The differential
defined here uses $N$-complexes. That puts right the wrong differential
presented in~\cite{rb:nonquad} in a 2-complex setting. Actually, as we
shall see it below, it is impossible to avoid $N$-complexes in
defining the differential, whereas the bimodule Koszul complex is
a 2-complex. Note that we shall keep the notation $s$ (instead of $N$) of~\cite{rb:nonquad}.
\end{abstract}

This short communication is a corrigendum to my article
``Koszulity for nonquadratic algebras''~\cite{rb:nonquad}.
In Section 5 of~\cite{rb:nonquad}, just below Lemma 5.5, it is not true
that the differentials $d'_{L}$ and $d'_{R}$ over $K_{L-R}$ commute if $s\geq 3$. The differential $d'$
has to be changed in the following manner. The new (true) differential
uses $s$-complexes. For $s$-complexes, see~\cite{bdvw:homog} and
references therein. 

Introduce
$(\overline{K}_{L}, \delta_{L})$ by $\overline{K}_{L,n}=A \otimes
J_{n}$, $n\geq 0$, and the $A$-linear map $\delta_{L} :
\overline{K}_{L,n} \rightarrow \overline{K}_{L,n-1}$ is defined by
the natural inclusion $J_{n} \hookrightarrow A\otimes J_{n-1}$. We
have $(\delta_{L})^{s}=0$, so $(\overline{K}_{L}, \delta_{L})$ is a
$s$-complex. Define analogously the $s$-complex $(\overline{K}_{R},
\delta_{R})$. Then $\overline{K}_{L-R}=\overline{K}_{L}\otimes A =
A\otimes \overline{K}_{R}$ is a bimodule $s$-complex for $\delta_{L}' =
\delta_{L}\otimes 1_{A}$ (respectively, for $\delta'_{R} =
1_{A} \otimes \delta_{R}$). But now $\delta'_{L}$ and $\delta'_{R}$
commute. The new differential $d'$ over $K_{L-R}$ is
defined by 
$$d'_{i}= \delta'_{L} - \delta'_{R},$$
if $i$ is odd, and
$$d'_{i}= \delta_{L}^{'s-1}+ \delta_{L}^{'s-2} \delta'_{R} + \ldots + 
\delta'_{L} \delta_{R}^{'s-2} + \delta_{R}^{'s-1},$$
if $i$ is even. So we obtain a pure projective complex in the
category $\mathcal{C}$, called the \emph{bimodule Koszul complex} of $A$. The
exact sequence (5.6) still holds, as well Theorem 5.6 with the same
proof.

Formulas (5.12) and (5.13) have to be replaced accordingly. Assuming
$i$ odd $\geq 3$, we have
$$\tilde{d}_{i}(a\otimes vwv') =av\otimes
w v' - \ v'a\otimes v w,$$
where $v$ and $v'$ are in $V$, and $w$ is in $V^{(js-1)}$. Assuming
$i$ even, we have
$$\tilde{d}_{i}(a\otimes v_{1}\ldots v_{js})= av_{1}\ldots
v_{s-1} \otimes v_{s} \ldots v_{js} + v_{js}av_{1}\ldots
v_{s-2}\otimes v_{s-1} \ldots v_{js-1}$$
$$ + \cdots + v_{js-s+2}\ldots
v_{js}a \otimes v_{1}\ldots v_{js-s+1},$$
where $v_{1}, \ldots , v_{js}$ are in $V$.

While (5.16) still holds, (5.15) is replaced by a sum of $s$
appropriate sums, each of those is indexed by $Q\subset P,\, |Q|=s-1$. In particular,
$s[P]$ takes place of $2[P]$ in (5.23). All the remainder is
unchanged.


\begin{thebibliography}{99}


\bibitem{rb:nonquad} R. Berger, Koszulity for nonquadratic algebras,
  \emph{J. Algebra}  \textbf{239} (2001),
705-734.
\bibitem{bdvw:homog} R. Berger, M. Dubois-Violette, M. Wambst,
  Homogeneous Algebras \emph{J. Algebra}, to appear. 

\end{thebibliography}
\end{document}